\documentclass[12pt]{amsart}
\usepackage{amsfonts, amssymb,amsmath, amsthm, amsxtra, latexsym, mathrsfs, amscd}
\usepackage{amsmath, amssymb, amsthm, times, float}
\usepackage{graphics}

\usepackage{amssymb,hyperref}

\usepackage[latin1]{inputenc} %

\newtheorem{theo+}{Theorem}[section]
\newtheorem{prop+}[theo+]{Proposition}
\newtheorem{coro+}[theo+]{Corollary}
\newtheorem{lemm+} [theo+]{Lemma}
\newtheorem{deep+}  [theo+]  {Deep Result}
\newtheorem{fact+}  [theo+]  {Fact}
\theoremstyle{definition}
\newtheorem{exam+}  [theo+]  {Example}
\newtheorem{rema+}  [theo+]  {Remark}
\newtheorem{defi+}  [theo+]  {Definition}
\newtheorem{xca+}[theo+]{Exercise}
\newenvironment{theorem}{\begin{theo+}}{\end{theo+}}

\newenvironment{remark}{\begin{rema+}}{\end{rema+}}

\numberwithin{equation}{section}

\usepackage{color}

\def\draft{\centerline{(Draft {\the \day}/{\the\month} \the \year.)}}
\bigskip

\def\refn#1.#2{\expandafter\def\csname#1\endcsname{[#2]}}
\def\refnr#1.{\csname#1\endcsname}

\def\a{\alpha}

\def\Claminv2{|C(\Lambda)|^{-2}}

\def\de{d\varepsilon}

\def\Aa2D{A^{\a,2}(D)}
\def\bAa2D{\overline{A^{\a,2}(D)}}
\def\Ab2D{A^{\beta,2}(D)}
\def\bAb2D{\overline{A^{\beta,2}(D)}}

\def\Norm#1_#2{\Vert#1\Vert_{#2}}

\def\phipl12{\phi_{p_{l_1}, p_{l_2}}}
\def\phip01{\phi_{p_{0}, p_{0}}}

\def\a{\alpha}

\def\Claminv2{|C(\Lambda)|^{-2}}

\def\Sig{\Sigma}

\def\de{d\varepsilon}

\def\Aa2D{A^{\a,2}(D)}
\def\bAa2D{\overline{A^{\a,2}(D)}}
\def\Ab2D{A^{\beta,2}(D)}
\def\bAb2D{\overline{A^{\beta,2}(D)}}

\def\phipl12{\phi_{p_{l_1}, p_{l_2}}}
\def\phip01{\phi_{p_{0}, p_{0}}}

\def\bc{\mathbb C}

\def\br{\mathbb R}

\def\alg/{algebra}

\def\Alg/{Algebra}

\def\alt/{alternative} 
\def\anal/{analytic}
\def\analfunc/{\anal/\ \func/}
\def\Ans/{\it Answer. \normal}
\def\ass/{associative}
\def\nass/{non-\ass/}
\def\autom/{automorphism}
\def\homom/{homomorphism}
\def\isom/{isomorphism}
\def\bdd/{bounded}
\def\Bdd/{Bounded}
\def\bddsymdom/{bounded \sym/ \dom/}
\def\Cartdom/{Cartan \dom/}
\def\bdry/{boundary}
\def\bsd/{\bdd/ \symdom/}
\def\bv/{boundary value}
\def\cf/{{\it cf}\.}
\def\Cf/{{\it Cf}\.}
\def\charr/{character}
\def\coeff/{coefficient}
\def\comm/{commutative}
\def\cpct/{compact}
\def\compl/{complex}
\def\comp/{complex}
\def\Comp/{Complex}
\def\conf/{conformal}
\def\conj/{conjugate}
\def\conn/{connect}
\def\cont/{continuous}
\def\conv/{converge} 
\def\convc/{convergence}
\def\convt/{convergent}
\def\convx/{convex}
\def\coord/{coordinate}
\def\lcoord/{local coordinate}
\def\Corr/{Corresponding}
\def\corr/{corresponding}
\def\corrd/{correspond}
\def\cov/{covariant}
\def\decomp/{decomposition}
\def\deco/{decompose}
\def\diff/{different} 
\def\Diff/{Different} 
\def\dimn/{dimension} 
\def\distr/{distribution} 
\def\div/{diverge} 
\def\dom/{domain}
\def\eg/{\hbox{\it e.g}\.}
\def\eigenf/{eigen\-\func/}
\def\eigensp/{eigen\-space}
\def\eigenv/{eigen\-value}
\def\eq/{equation}
\def\equa/{equation}
\def\de/{\diff/ial \equa/}
\def\do/{\diff/ial operator}
\def\ode/{ordinary \de/}
\def\pde/{partial \de/}
\def\pdo/{partial \diff/ial operator}
\def\psdo/{pseudo \diff/ial operator}
\def\fin/{finite}
\def\Ex/{\it Example.\ \normal}
\def\Exnr#1/{\it Example #1.\ \normal}
\def\foll/{follow}
\def\follg/{following}
\def\Follg/{Following}
\def\func/{function}
\def\Func/{Function}
\def\Fonc/{Fonc\-tion}
\def\fonc/{fonc\-tion}
\def\Funk/{Funk\-tion}
\def\funk/{Funk\-tion}
\def\gen/{general}
\def\har/{harmonic}
\def\Hint/{\it Hint. \normal}
\def\hist/{historic}
\def\histcl/{historical}
\def\hol/{holo\-morphic}
\def\homog/{ho\-mo\-ge\-ne\-ous}
\def\hyp/{hyper\-bolic}
\def\hyperg/{hyper\-geometric}
\def\ie/{\hbox{\it i.e.}}
\def\iff/{if and only if}
\def\ineq/{inequality}
\def\infra/{{\it inf\-ra}}
\def\ultra/{{\it ult\-ra}}
\def\Inpart/{In particular}
\def\inpart/{in particular}
\def\instof/{instead of}
\def\interps/{interpolation space}
\def\interp/{interpolation}
\def\Interp/{Interpolation}
\def\interpr/{Interpretation}
\def\Intr/{Introduction}
\def\intv/{interval}
\def\inv/{invariant}
\def\invc/{invariance}
\def\Iowords/{In other words}
\def\iowords/{in other words}
\def\ipr/{inner product}
\def\irred/{irreducible}
\def\lb/{line bundle}
\def\lin/{linear}
\def\lhs/{left hand side}
\def\rhs/{right hand side}
\def\loc/{local}
\def\math/{mathematic}
\def\mathcn/{\math/ian}
\def\manif/{manifold}
\def\meas/{measure}
\def\measl/{measurable}
\def\mero/{mero\-morphic}
\def\mon/{monomial}
\def\monog/{monogenic}
\def\mult/{multiple}
\def\multy/{multiply}
\def\multn/{multiplication}
\def\nas/{necessary and sufficient}
\def\nbd/{neighborhood}
\def\neg/{negative}
\def\nondeg/{nondegenerate}
\def\Oohand/{On the other hand}
\def\oohand/{on the other hand}
\def\Oonhand/{On the one hand}
\def\oonhand/{on the one hand}
\def\oper/{operator}
\def\orth/{ortho\-gonal}
\def\orthon/{ortho\-normal}
\def\otoh/{on the other hand}
\def\quat/{quaternion}
\def\pp/{\hbox{a. e.}}
\def\psh/{plurisubharmonic}
\def\pol/{polynomial}
\def\pot/{potential}
\def\pos/{positive}
\def\princ/{principle}
\def\prob/{probability}
\def\proj/{projective}
\def\projn/{projection}
\def\Proof/{\it Proof:\normal}
\def\Rem/{\it Remark\normal}
\def\Remnr#1/{\it Remark\ \normal #1. }
\def\rep/{representation}
\def\reps/{representations}
\def\meta/{metaplectic representation}
\def\repr/{reproducing}
\def\reprker/{reproducing kernel}
\def\resp/{respective} 
\def\resply/{respectively}
\def\restr/{restriction}
\def\sa/{self-adjoint}
\def\st/{such that}

\def\sol/{solution}
\def\ru/{space}
\def\sph/{spherical}
\def\ssp/{sub\ru/}
\def\sym/{symmetric}
\def\Sym/{Symmetric}
\def\symb/{symbol}
\def\symbc/{symbolic}
\def\symdom/{\sym/ domain}
\def\symp/{symplectic}
\def\Theor#1/{\fet Theorem #1.\ \normal}
\def\Lem#1/{\fet Lemma #1.\ \normal}
\def\Lemma/{\fet Lemma.\ \normal}
\def\topl/{topology}
\def\topll/{topological}
\def\transf/{transform}
\def\transl/{translation}
\def\transfn/{transformation}
\def\transv/{transvectant}
\def\trig/{trigonometric}
\def\tril/{trilinear}
\def\trilf/{trilinear form}
\def\uhp/{upper halfplane}
\def\uhs/{upper halfspace}
\def\vb/{vector bundle}
\def\vf/{vector field}
\def\vsp/{vector space}
\def\wrt/{with respect to}
\def\Wlog/{Without loss of generality}
\def\a{\alpha}

\def\Sig{\Sigma}

\def\Ab/{Abel}
\def\Ban/{Banach}
\def\Bansp/{\Ban/ space}
\def\Belt/{Bel\-tra\-mi}
\def\Berg/{Berg\-man}
\def\Bern/{Ber\-nou\-lli}
\def\Berz/{Berezin}
\def\Bess/{Bessel}
\def\Cart/{Car\-tan}
\def\Cay/{Cay\-ley}
\def\CG/{Clebsch-Gordan}
\def\Cl/{Clifford}
\def\CR/{Cauchy-Rie\-mann}
\def\Dir/{Dirichlet}
\def\Eucl/{Euclide}
\def\Eucln/{Euclidean}
\def\F/{Fourier}
\def\Hank/{Hankel}
\def\Hankf/{\Hank/ form}
\def\Herm/{Hermite}
\def\Hilb/{Hilbert}
\def\Hilbs/{Hilbert space}
\def\Hilbsp/{Hilbert space}
\def\HS/{Hilbert-Schmidt}
\def\Lag/{La\-grange}
\def\Lap/{La\-place}
\def\LapBelt/{\Lap/-\Belt/}
\def\Leb/{Lebesgue}
\def\Marc/{Mar\-cin\-kie\-wicz}
\def\Moeb/{Moebius}
\def\Moebt/{Moebius transformation}
\def\Moebtransfn/{Moebius transformation}
\def\Pla/{Plan\-che\-rel}
\def\Poin/{Poin\-car\'e}
\def\Riem/{Rie\-mann}
\def\Riemn/{\Riem/ian}
\def\psRiemn/{pseudo-\Riem/ian}
\def\Riems/{Rie\-mann surface}
\def\Schroe/{Schr\"odinger}
\def\Weier/{Weier\-strass}
\def\anal/{analytic}
\def\bsd/{bounded symmetric domain  }
\def\bdd/{bounded}
\def\calc/{calculation}\def\conj{conjugate}
\def\calci/{calculating}\def\eg{e.g.}
\def\conj/{conjugate}
\def\deco/{decomposition}
\def\eg/{e.g.}
\def\fct/{function}
\def\gp/{group}
\def\hw/{highest weight}
\def\hwv/{highest weight vector}
\def\hwvs/{highest weight vectors}
\def\lw/{lowest weight}
\def\lwv/{lowest weight vector}
\def\lwvs/{lowest weight vectors}
\def\hds/{holomorphic discrete series}
\def\iff/{if and only if}
\def\inv/{invariant}
\def\irrde/{irreducible decomposition}
\def\meas/{measure}
\def\transf/{transform}
\def\rep/{representation}
\def\resp/{respectively}
\def\inters/{intertwines}
\def\interg/{intertwining}
\def\meta/{metaplectic representation}
\def\qu/{quaternion}
\def\rep/{representation}
\def\symdom/{ symmetric domain}
\def\st/{such that}
\def\shd/{subhead}
\def\transf/{transform}
\def\wrt/{with respect to}

\def\ra{\rightarrow}

\def\Norm#1#2#3{\Vert#1\Vert^{#3}_{{#2}}}

\begin{document}

\def\abstractname{Abstract}
\def\chrefname{References}

\title[Hitchin component
]{K\"ahler metric on the space of convex real projective structures on surface
}
\author{Inkang Kim and  Genkai Zhang}

\address{School of Mathematics,
KIAS, Heogiro 85, Dongdaemun-gu
Seoul, 130-722, Republic of Korea.
\text{Email: inkang@kias.re.kr}
}
\address{Mathematical Sciences, Chalmers University of Technology and
Mathematical Sciences, G\"oteborg University, SE-412 96 G\"oteborg, Sweden.
\text{Email: genkai@chalmers.se}
}
\footnotetext[3]{Research partially supported by
STINT-NRF grant (2011-0031291). Research by G. Zhang is supported partially
 by the Swedish
Science Council (VR). I. Kim gratefully acknowledges the partial support
of grant  (NRF-2014R1A2A2A01005574) and a warm support of
Chalmers University of Technology during his stay.
}
\begin{abstract}
We prove  that the space of convex real projective structures on
a surface of genus $g\ge 2$ admits a mapping class group invariant
  K\"ahler
metric where Teichm\"uller space with  Weil-Petersson metric  is a totally geodesic
complex submanifold.
\end{abstract}

\maketitle

\baselineskip 1.35pc

\section{Introduction}
Recently, the character variety $\chi(\pi_1(M), G)$ of representations of $\pi_1(M)$ in   a real algebraic group $G$, has drawn many attentions from different branches of mathematics. The $G$-character variety  $\chi(\pi_1(M), G)$ is the geometric quotient
of $Hom(\pi_1(M), G)$ by inner automorphisms of $G$. Often, some components of the character variety correspond to some geometric structures on $M$.
Hitchin \cite{Hitchin} introduced  Hitchin component in the
character variety of a closed surface group in $PSL(n+1,\mathbb R)$
generalizing Teichm\"uller space in $PSL(2,\mathbb R)$. More precisely, he showed that for the adjoint group $G$ of the
split real form of a complex simple Lie group $G^c$, the quotient by the conjugation action of $G$ of the set of
homomorphisms, from the fundamental group $\Gamma$ of a closed surface $S$ of
genus $g\ge 2$ to $G$, which acts completely reducibly
on the Lie algebra of $G$, has a connected component homeomorphic to  Euclidean space of dimension $(2g-2)\mathrm{dim} G$.
His method is the use of Higgs bundle theory developed by himself, K. Corlette, S. Donaldson, C. Simpson and many others \cite{Hit,Sim}.
A homomorphism from $\Gamma$ to $G^c$ defines a flat principal $G^c$-bundle. Given a complex structure on $S$, denoted by $\Sigma$, a
theorem of Corlette and Donaldson associates a natural $G^u$-connection
$A$, where $G^u$ is the maximal compact subgroup of $G^c$, and a Higgs
field $\Psi\in H^0(\Sigma, \mathrm{ad} P\otimes\mathcal K)$ which
satisfy the
equation $F_A+[\Psi,\Psi^*]=0$. Here $\mathcal K$ is the canonical
line bundle over $\Sigma$, $P$ is a principal $G^u$-bundle and $\mathrm{ad} P$ is the Lie algebra bundle associated to the adjoint representation of $G^u$. Solutions to these equations provide a holomorphic parametrization of the equivalence classes of
homomorphisms from $\Gamma$ to $G^c$. For appropriately
chosen $\Psi$ the solutions are stable under the complex conjugation
in $G^c$ and reduces to $G$-connection corresponding to elements
in the Hitchin component.

For the real linear group $SL(n+1,\br)$, Hitchin  showed that the Hitchin component is homeomorphic to $\bigoplus_{j=1}^n H^0(\Sigma,\mathcal K^{j+1})$.
We need to mention that this homeomorphism depends on a priori fixed
 complex structure on $S$, and hence it   is not
mapping class group equivariant.
After Hitchin's work, many people pursued to clarify this component in many different ways. Notably Labourie \cite{La1} introduced a
notion of Anosov representations and proved that Hitchin representations are exactly Anosov representations in $SL(n+1,\br)$.
In \cite{La2}, he also suggested a mapping class group equivariant parametrization using  Hitchin map and an adaptation of an
energy functional over Teichm\"uller space.  We will review his interpretation of Hitchin map in Section \ref{Hitchin}.

It has been conjectured for a long time that the Hitchin component admits a mapping class group invariant K\"ahler metric.
There have been many evidences for this, see \cite{G-crelle, Li, BCLS}.
In the last section, we prove the existence of a mapping class group invariant K\"ahler metric on the Hitchin component for $n=2$.
\begin{theo+}
The Hitchin component of the character variety
$\chi(\pi_1(M),
SL(3,\br))$
can be equipped with a mapping class group
invariant  K\"ahler metric where $M$ is a closed surface of genus $\geq 2$.
Furthermore  Teichm\"uller space equipped with
the Weil-Petersson metric is a totally geodesic complex submanifold.
\end{theo+}
This K\"ahler metric is constructed using
certain $L^2$-metric.  Intuitively, we glue Weil-Petersson metric on the base and $L^2$-metric along vertical fibers using Griffith negativity. Indeed, we need the dual of a holomorphic vector bundle over Teichm\"uller space whose fibres are cubic holomorphic forms.
The geometric properties of this metric such as various curvatures, geodesics will be explored in a near future. We hope that this new natural K\"ahler metric will help us to better understand the moduli space of real projective structures and the Teichm\"uller space  as a byproduct.

This particular Hitchin component has been intensively studied by many people.
 Choi-Goldman showed that the corresponding geometric structure is
the convex real projective structure \cite{CG} and the bundle structure is verified by Labourie \cite{La} and Loftin \cite{Loftin} independently using
Monge-\`Ampere equations relying on the seminal work of Cheng-Yau \cite{CY}. The symplectic structure
on the Hitchin component has
been studied by Goldman \cite{G-crelle}.
 Recently Li \cite{Li}
 constructed a mapping class group invariant metric using  explicit constructions
over $\mathfrak{sl}(3,\br)$-bundles and Cheng-Yau metric over the
cone.  It would be interesting to compare this construction with ours.
Bridgeman-Canary-Labourie-Sambarino \cite{BCLS} constructed a pressure
metric on Hitchin component of $SL(n+1, \mathbb R)$
 for every $n$ using different method.

{\bf Acknowledgement}
We would like to thank Bo Berndtsson for a few helpful discussions
on complex vector bundles, in particular for  K\"a{}hler property
based on Griffith's positivity. We are grateful
for the anonymous referee for valuable suggestions
and for the careful reading of an earlier version
of this paper.

\section{Convex projective structures on a manifold $M$}
A flat projective structure on an $n$-dimensional manifold $M$ is a $(\mathbb{RP}^n, PSL(n+1,\br))$-structure, i.e., there exists a maximal atlas on $M$ whose transition maps are restrictions to open sets in $\mathbb{RP}^n$ of elements in $PSL(n+1,\br)$. Then there exist a natural holonomy map $\rho:\pi_1(M)\rightarrow PSL(n+1,\br)$ and a developing map from the universal cover $\tilde M$, $f:\tilde M\rightarrow \mathbb{RP}^n$ such that
$$\forall x\in\tilde M,\ \forall \gamma\in\pi_1(M), \ f(\gamma x)=\rho(\gamma)f(x).$$
We will consider projective structures deformed from  hyperbolic structures,  and all holonomy representations will lift to $SL(n+1,\br)$.
An $\mathbb{RP}^n$-structure is convex if the developing map is a homeomorphism onto a convex domain in $\mathbb{RP}^n$.  It is {\it properly convex} if the domain is included in a compact convex set of an affine chart, {\it strictly convex} if the convex set is strictly convex.

When $M=S$ is a closed Riemann surface of genus at least 2, a huge amount of literature for the set of marked strictly convex real projective structures on $S$ exist  concerning its parametrization \cite{Gol},
its identification with Hitchin component  \cite{CG},   degeneration of the projective structures \cite{Kim1},  entropy of geodesic flow \cite{crampon}, the marked length rigidity \cite{Kim}, its Zariski tangent space at Fuchsian locus of the character variety in $SL(n,\br)$ \cite{KZ1} and many more. For a recent generalization to  finite volume convex real projective structures, see \cite{Benoist, CLT}.

In this paper, we utilize the holomorphic vector bundle structure of the space of the marked strictly convex real projective structures on a closed surface of genus at least 2, \cite{La, Loftin},
where Monge-\`Ampere equation type argument is used. The method is initiated by  Cheng and Yau \cite{CY, CY1}.

\section{Hitchin map and bundle structure}
\label{Hitchin}
In this section we collect known results to introduce the Hitchin map. See \cite{La, La2} for details.
If $\rho$ is a reductive representation from
$\pi_1(S)$ to a semisimple Lie group $G$,  Corlette proved
\cite{Co} the following claim: There exists a unique, up to $G$, $\rho$-equivariant harmonic map
$f_{\rho,J}:\tilde S \ra G/K$ where $G/K$ is
the symmetric space of $G$ and $J$ is a complex structure on $S$. If $f$ is a $\rho$-equivariant map from the universal
cover $\tilde S$ to $G/K$, then one can define the energy $E(J,f)$ of $f$ with respect to a complex structure $J$ on $S$. Hence given $\rho$, one can define an energy functional $e_\rho(J)$ on Teichm\"uller space by the infimum
of energy of $\rho$-equivariant functions with respect to a complex structure $J$. This harmonic map minimizes the energy, i.e. $E(J,f_{\rho,J})=e_\rho(J)$. The minimum area of
$\rho$ is defined to be $\inf_{J}e_\rho(J)$. Then it is known \cite{SU, SY} that a harmonic map realizing the minimum area of $\rho$ is conformal.

For a map $f$ from $S$ to a Riemannian manifold $(M,g)$, $Tf$ can be viewed as a 1-form on $S$ with values in $f^*TM$. Let $T_\bc f(u)=Tf(u)-iTf(Ju)$ be a complexified tangent map. Then f is harmonic if and only if $T_\bc f$ is holomorphic. Furthermore $g_\bc(T_\bc f, T_\bc f)=0$ if and only if $f$ is minimal.

Every $G$-invariant symmetric multilinear form $P$ on $\mathfrak g$ gives rise to a parallel polynomial function $P$, with the same notation, on $G/K$.
Hence for any complex structure $J$ on $S$, and for every symmetric  $G$-invariant multilinear form $P$ of degree $k$ on $\mathfrak g$, any reductive representation $\rho$ gives rise to an element in $\mathcal Q(k,J)$ by
$P(T_\bc f,\cdots,T_\bc f)$ where $f$ is a $\rho$-equivariant harmonic map. Here $\mathcal Q(k,J)=H^0((S,J),\mathcal K^k)$ is the space of holomorphic $k$-differentials. Denote this map by $F_{P,J}(\rho)=P(T_\bc f,\cdots,T_\bc f)$.

For  $G=SL(n,\br)$, we can use the symmetric polynomial $P_k$ of degree $k$. Then $F_{P_2}$ is a metric on $SL(n,\br)/SO(n)$.  Set
$$\Psi_J=\bigoplus_{k=2}^{k=n} F_{P_k,J}.$$

Hitchin proved that the map $\Psi_J$ is a homeomorphism from the Hitchin component to $\mathcal Q(2,J)\oplus \cdots \oplus Q(n,J)$. Set $\epsilon^{(n)}$ to be the bundle over Teichm\"uller space with fibres $\mathcal Q(k,J),\ k\geq 3$. Labourie introduced the Hitchin map from $\epsilon^{(n)} $ to the Hitchin component
$$H(J,\omega)= \Psi_J^{-1}(\omega)$$ for $\omega\in \bigoplus_{k=3}^{n} Q(3,J)\oplus \cdots \oplus Q(n,J)$.
He showed that
\begin{theo+}
The Hitchin map is surjective.
\end{theo+}

For $n=3$, this Hitchin map is injective also.
Let  $\mathcal T$ be the Teichm\"uller space
of complex structures $\Sig_t$ on the surface $S$
with the holomorphic tangent space given by
$H^{(0, 1)}(\Sig_t, \mathcal K^{-1})$ at each $t\in \mathcal T$.
Following \cite{La} we consider
the space
$$
\mathcal V=\{
(v, t);
v\in V_{t}:=
{H^0(\mathcal K_t^3)},
\,  t\in \mathcal T\}.
$$
Labourie \cite{La} showed that
\begin{theorem}There exists a mapping class group equivariant diffeomorphism between
the moduli space of convex structures on $S$ and the moduli space of pairs $(J,Q)$ where $J$ is  a complex
structure on $S$ and $Q$ is a cubic holomorphic differential on $S$
with respect to $J$,
i.e. diffeomorphic to $\mathcal V$.
\end{theorem}
Hence the moduli space of convex projective structures can be treated
as
the  holomorphic vector bundle $\mathcal V$
over $ \mathcal T$.

The vector bundle $\mathcal V$ (sometimes called Hodge bundle)
can be treated as in \cite{Bob}.
Consider the tautological bundle \cite{Ahlfors-1961}
$$
\hat T
=\{(t, w); t\in
{\mathcal T},
 w\in \Sig_t\}$$
over the Techm\"u{}ller space ${\mathcal T}$
 equipped with the
canonical complex structure; see
\cite[\S5]{Ahlfors-1961}.
Note that $\hat T
$
is  then a K\"a{}hler manifold
 with the K\"a{}hler metric
being locally the product metric.

The Hodge vector bundle $
\mathcal V
$
is then a complex vector bundle and in particular it is
a complex manifold, see e.g. \cite{Bob}.
We let $\mathcal V^\ast$ be the dual bundle of $\mathcal V$.
This can be realized as
$$
\mathcal V^\ast=\{
(v, t);
v\in V_{t}:=
{H^{0, 1}(\mathcal K_t^{-2})},
\,  t\in \mathcal T\}
$$
via the natural mapping class group invariant
paring $(f, g)=\int_{\Sig_t}g(f)$, $f\in
{H^{0, 1}(\mathcal K_t^{-2})},
g\in {H^{0}(\mathcal K_t^{3})}
= {H^{1}(\mathcal K_t^{2})}$.


\section{Various notions of curvature positivity}
We shall need some results on the positivity for the curvature
of the Hodge bundle above. Recall that generally
a complex  manifold $M$ equipped with
a Hermitian metric is said to have a {\it
  nonpositive bisectional curvature} if
$$R(X,Y,\bar X,\bar Y)\leq 0$$ for all $X,Y\in TM\otimes\mathbb C$
where $R$ is the curvature tensor extended complex linearly to
complexified bundle.
When we deal with the holomorphic vector bundles, there are similar notions of positivity (negativity).
Let $E$ be a holomorphic vector bundle over a  K\"ahler
manifold $M$ and $h$ a Hermitian metric on $E$. Let $\nabla$ be a
Chern connection which is compatible with the metric $h$ and complex
structure on $E$. If we write $\nabla=\mathcal D+
\bar \partial$, then its curvature $F$ is equal to $\mathcal D\bar
\partial-\bar \partial\mathcal D$ when acting
on local holomorphic sections.
 More concretely,
for any section $s$ and (complexified) vector fields $X,Y$,
$$F(X,Y)(s)=\nabla_X\nabla_Y s-\nabla_Y\nabla_X s- \nabla_{[X,Y]} s.$$
$F$  is of type $(1,1)$, real and satisfies
$$h(F(s_1), s_2)+h(s_1,F(s_2))=0$$ for any sections $s_1$ and $s_2$.

If $z^i$ are local holomorphic coordinates on $M$ and $e_\alpha$ is a local orthogonal  frame on $E$, then the curvature $F$ can be written by
$$\sqrt{-1} F=\sum c_{j\bar k}^{\bar{\alpha}\beta} dz^jd\bar z^k\otimes e^*_\alpha \otimes e_\beta,$$ where $\overline{ c_{j\bar k}^{\bar\alpha\beta}}=c_{kj}^{\bar\beta\alpha}$. This curvature gives rise to a Hermitian (sesqui-linear) form
 $\Theta$  on $TM\otimes E$, given locally by
$$\Theta=\sum  c_{j\bar k}^{\bar{\alpha}\beta} (dz^j\otimes e^*_\alpha)\otimes\overline{(dz^k\otimes e^*_\beta)}.$$

In  tensorial notation, let
 $e_\alpha$ be a
local holomorphic frame of $E$
and $e^{\alpha}$ the dual frame, then the curvature tensor $R\in \Gamma(M, \wedge^2T^*M\otimes E^*\otimes E)$ of $\nabla$ has the form
$$R=\frac{\sqrt{-1}}{2\pi}\sum R^\gamma_{i\bar j \alpha} dz^id\bar z^j\otimes e^\alpha\otimes e_\gamma$$
 where  $R^\gamma_{i\bar j \alpha}=h^{\gamma\bar\beta}R_{i\bar j\alpha\bar\beta}$ and
$$R_{i\bar j\alpha\bar\beta}=-\frac{\partial^2 h_{\alpha\bar\beta}}{\partial z^i\partial\bar z^j}+\sum h^{\gamma\bar\delta}\frac{\partial h_{\alpha\bar\delta}}{\partial z^i}\frac{\partial h_{\gamma\bar\beta}}{\partial\bar z^j}.$$

Then the Hermitian vector bundle $(E,h)$ is said to be
\begin{enumerate}
\item Griffith positive if for any nonzero vectors $u=\sum u^i\frac{\partial}{\partial z^i}$ and $v=\sum v^\alpha e_\alpha$,
$$\sum R_{i\bar j\alpha\bar\beta} u^i\bar u^j v^\alpha \bar v^\beta>0,$$ i.e.,
$$ \Theta(u\otimes v, u\otimes v)=h(F(u, \bar u)(v),v)>0$$ for any section $v\neq 0$ and non-zero holomorphic tangent vector field $u$. In other words, $\Theta$ is positive definite on nonzero simple tensors of the form $u\otimes v$.
\item Nakano positive if for any nonzero vector $u=\sum u^{i\alpha}\frac{\partial}{\partial z^i}\otimes e_\alpha$,
$$\sum  R_{i\bar j\alpha\bar\beta} u^{i\alpha}\bar u^{j\beta}>0,$$  i.e.,
the associated sesqui-linear form $\Theta$ is a positive definite Hermitian form.
\item dual Nakano positive if for any nonzero vector $u=\sum u^{i\alpha}\frac{\partial}{\partial z^i}\otimes e_\alpha$,
$$\sum  R_{i\bar j\alpha\bar\beta} u^{i\beta}\bar u^{j\alpha}>0.$$
\end{enumerate}
It is obvious from the definition that  Nakano positivity implies Griffith positivity.
Also $(E,h)$ is dual Nakano positive if and only if
$(E^*,h^*)$ is  Nakano negative. If $E$ is Griffith positive then its dual $E^*$ is Griffith negative.
We also remark that Griffith positivity implies that the bundle $E$ is ample.

In particular if $E$ is a line bundle, then the above notions of
positivity  all
agree. The curvature in this case is computed by
$\bar\partial\partial \log \Vert w\Vert^2$ for any local holomorphic
frame (i.e. non zero-section) $w$.
We will use these notions of positivity to prove the existence of K\"ahler metric on the Hitchin component of $\chi(\pi_1(S),SL(3,\mathbb R))$.

\section{K\"ahler property}

The following result is a corollary of a general
theorem of Berndtsson \cite{Bob} applied to the vector bundle
$\mathcal V$; when we replace the fiber of the bundle $\mathcal V$
by  the spaces $H^0(\mathcal K_t^2)$ viewed as the dual
space of the holomorphic tangent space of Teichm\"uller
space, it is
a classical result of Ahlfors \cite{Ahlfors-1961/62} that
the Teichm\"uller space has negative holomorphic
sectional curvature equipped with the Weil-Petersson metric.

\begin{theo+} The bundle $\mathcal V$ is
 Griffith positive.
\end{theo+}
\begin{proof} The tangent  line bundle $\mathcal K_t^{-1}$ on  each  Riemann surface $\Sigma_t$
is equipped with a unique K\"ahler-Einstein metric of negative curvature $-1$, in other words,
the canonical line bundle  $\mathcal K_t
$ is
K\"ahler-Einstein  of positive curvature $1$. Let $\mathcal L$ on $\hat T$
be the pull-back of the line bundle $\mathcal K_t$ to $\hat T$ under the projection $\hat T\rightarrow \mathcal T$.
 It follows from \cite[Theorem 5.5, Lemma 5.8]{Wolpert-invmath}
that the bundle  $\mathcal  L$ is
positive; see also  \cite[Main Theorem, Theorem 1]{Sch} for generalization.
Thus  $\mathcal  L^2$ is positive since $ \mathcal L$
is a line bundle.
In \cite{Bob}, it is proved that the Hermitian vector bundle over $\mathcal T$,
$$
H^0(\Sigma_t, L_0\otimes \mathcal K_{\hat T/\mathcal T})
 \mapsto t
$$
the fiber being the spaces of global sections,
endowed with $L^2$-metric, is Nakano positive,
where $L_0$ is any positive line bundle over the K\"ahler manifold $\hat T$.
In particular taking $L_0=\mathcal L^2$ we have
$$
H^0(\Sigma_t, \mathcal L^2\otimes \mathcal K_{\hat T/\mathcal T})
=H^0(\Sigma_t, \mathcal K_t^2\otimes \mathcal K_{\hat T/\mathcal T})
=H^0(\Sigma_t, \mathcal K_t^3)$$
we see that
the bundle with fiber $H^0(\mathcal K_t^3)$ over $\mathcal T$,
i.e., the  bundle $\mathcal V$ is Nakano positive, hence Griffith
positive.
\end{proof}

\begin{rema+} The Nakano positivity of the bundle
$H^0(\Sigma_t, \mathcal K_t^3)$ can also be directly
proved by using the result of Berndtsson in \cite{Bo-mz}.
It is proved there whenever the metric on $L_0$ is positive
fiberwise then the Nakano positivity still holds for
$H^0(\Sigma_t, L_0\otimes \mathcal K_{\hat T/\mathcal T})$. In our
case the fiber metric  is dual to the hyperbolic metric
on the Riemann surfaces and has constant  positive  curvature. See also \cite{LSY, LSY1} for related works.
\end{rema+}

\begin{coro+} The bundle $\mathcal V^*$ is
 Griffith negative.
\end{coro+}

The first statement below can be proved for general bundles with
the Griffith negativity, and here we are only interested in
the special case of  $\mathcal V^*$.

\begin{theo+}\label{main}
The bundle $\mathcal V^*$
is
a K\"ahler manifold.
In particular the Hitchin component of the character variety
$\chi(\pi_1(S),
SL(3,\br))$
 has a mapping class group invariant K\"a{}hler metric.
\end{theo+}

\begin{proof} Let  $w_0\ne 0$ be a fixed point in $
\mathcal V^*$ with $z_0=\pi(w_0)$, where $\pi: \mathcal V^*
\to \mathcal T$ is the defining projection.
We choose a local trivializing holomorphic
frame $\{e_{\alpha}=e_{\alpha}(z)\}$ in a coordinate neighborhood $U$
of $z_0$, and write, with some abuse of notation, $z=(z^1, \cdots, z^n)$
as the  coordinate of $U$.  The local holomorphic coordinates near
$w_0$
 will
be $w=\sum_\alpha x^\alpha e_\alpha(z)\to (z^i,x^\alpha)$, and the holomorphic tangent
vectors are $T=(u, v)$ with $u=\sum u_i \frac{\partial}{\partial z^i}$
and  $v=\sum v_\alpha \frac{\partial}{\partial x^\alpha}$.

We let $\psi$ be a local K\"a{}hler potential for $\mathcal T$ near
$z_0$, $\mathcal T$ being equipped with the Weil-Petersson metric.
Thus $\pi^*\psi(w)$ is defined in a neighborhood of $w_0$. We let
$$
\phi(w)
=
\Vert w\Vert^2
+\psi(\pi(w))
=\Vert w\Vert^2
+\pi^\ast\psi(w)
$$
to be defined in a neighborhood of $w_0$.

 Let $\nabla =\mathcal D +\bar\partial$
be the Chern connection acting on sections of $\mathcal V^\ast$ and $\mathcal D e_\alpha=\sum\theta^\beta_\alpha e_\beta$ where $\theta^\beta_\alpha$ is type $(1,0)$. We fix
$T=(u,v)\ne 0$ at $w_0$ and perform
the differentiation
$\bar \partial_T\partial_T \phi$ with $\partial_T=\partial_u +\partial_v$.
The vectors $w=\sum_\alpha x^\alpha e_\alpha(z)$ will be viewed as holomorphic sections for fixed
$x_\alpha$. The curvature $R(u, \bar u)$  of Chern connection $\nabla$ is
$$R(u,\bar u)w=\nabla_u\nabla_{\bar u}w-\nabla_{\bar u}\nabla_u w-\nabla_{[u,\bar u]} w=
\nabla_u \bar\partial_u w-\bar\partial_u \mathcal D_u w= -\bar\partial_u \mathcal D_u w$$ using
$\bar\partial_u w=0$ and $[u,\bar u]=0$.
 We have
$$
\partial_u \Vert w\Vert^2=(\mathcal D_u w, w)+(w,\bar\partial_u w)
=(\mathcal D_u w, w), $$$$
\bar\partial_u\partial_u \Vert w\Vert^2 =\bar\partial_u(\mathcal D_u w, w)=(\bar\partial_u\mathcal D_u w, w)+(\mathcal D_u w,\mathcal D_u w)
=-(R(u, \bar u) w, w)
+ (\mathcal D_u w, \mathcal D_u w).
$$
Since $$\mathcal D_u w=\sum D_u x^\alpha e_\alpha+ \sum x^\alpha \mathcal D_u e_\alpha=\sum x^\alpha \theta^\beta_\alpha(u) e_\beta(z),$$
$\theta^\beta_\alpha(u)e_\beta(z)$ depends smoothly only on $z$, and hence
$\bar\partial_v \mathcal D_u w=0$. Also $\partial_v w=v$, hence
$$
\bar\partial_v\partial_u \Vert w\Vert^2
=(\mathcal D_u w, v),
 \quad \bar \partial_u \partial_v \Vert w\Vert^2
=(v, \mathcal D_u w)
$$
and
$$
\bar\partial_v\partial_v \Vert w\Vert^2
=(v, v).
$$
Thus
\begin{eqnarray*}
\bar \partial_T\partial_T
\Vert w\Vert^2
=
-(R(u, \bar u) w, w)
+ (\mathcal D_u w, \mathcal D_u w)
+(\mathcal D_u w, v) +
(v, \mathcal D_u w)+
(v, v)
\\
=-(R(u, \bar u) w, w)
+ (\mathcal D_u w +v, \mathcal D_u w+v)
\ge 0.
\end{eqnarray*}
Here we use Griffith negativity of $\mathcal V^*$ to have $(R(u, \bar u) w, w)<0$.
We have also
$$
\bar \partial_T\partial_T
\pi^\ast\psi
=\bar \partial_u \partial_u \psi,
$$
which is positive definite in $u$.
Thus
$$
\bar \partial_T\partial_T \phi
=
-(R(u, \bar u) w, w)
+ (\mathcal D_u w +v, \mathcal D_u w+v)
+\bar \partial_u \partial_u \psi =I + II+ III,
$$
a sum of two semi-positive Hermitian forms.
We prove that  it is positive at $w=w_0$.
Suppose the quadratic form
 vanishes.
 Then $I=II= III=0$.
If $w_0\ne 0$, then $I=0$ implies that   $u=0$ by Griffith negativity, and the second term
is $0=II=(v, v)$, which implies further that  $v=0$.
Suppose $w_0=0$, i.e.,   in the zero section. Then $III=0$
implies that $u=0$, which in turn implies
$0=II= (v, v)$ and $v=0$. In either cases, $T=(0,0)$, a contradiction to the choice of $T=(u,v)\neq 0$.
\end{proof}

It follows from the definition that the map $w\to -w$ is an isometry
and its fixed point set is the space of zero sections, namely the
Teichm\"uller space identified as a submanifold. We have
thus

\begin{coro+} Let Teichm\"uller space $\mathcal T$
be equipped with the mapping class group invariant Weil-Petersson metric and
the Hitchin component
$\chi_H(\pi_1(S),
SL(3,\br))$
  of the character variety
$\chi(\pi_1(S),
SL(3,\br))$
be equipped with the mapping class group invariant K\"ahler metric as in Theorem
\ref{main}.  Then $\mathcal T$
is a totally geodesic submanifold of $\chi_H(\pi_1(S),
SL(3,\br))$.
\end{coro+}

Recenlty Labourie \cite{L} generalized this theorem to the Hitchin components associated to all real split simple Lie groups of rank 2.

\begin{remark}
If we consider  the bundle $\mathcal W$  over Teichm\"uller space whose fiber over $t$ equal to $\sum_{j\geq 2} H^0(\Sigma_t, \mathcal K_t^j)$,  then it is
Griffith positive and the same method applies to show that the total space $\mathcal W^*$ has a mapping class group invariant K\"ahler metric. For example, the space of the marked complex projective structures on a closed surface of genus at least 2 is a holomorphic vector bundle over Teichm\"uller space with fibers being the space of holomorphic quadratic differentials.
Hence this space has a natural mapping class group invariant K\"ahler metric.
\end{remark}

\end{document}